\documentclass[a4paper,12pt]{amsart}
\usepackage{amssymb}
\usepackage{ifthen}
\usepackage{graphicx}
\usepackage{float}
\usepackage{caption}
\usepackage{subcaption}
\usepackage{cite}
\usepackage{amsfonts}
\usepackage{amscd}
\usepackage{amsxtra}
\usepackage{color}

\newtheorem{thm}{Theorem}
\newtheorem{cor}{Corollary}
\newtheorem{lem}{Lemma}
\newtheorem{rem}{Remark}

\newtheorem{conj}{Conjecture}
\newtheorem{prob}{Problem}

\theoremstyle{definition}
\newtheorem{defn}{Definition}[section]
\newtheorem{example}{Example}

\newenvironment{pf}[1][]{%
 \vskip 1mm
 \noindent
 \ifthenelse{\equal{#1}{}}%
  {{\slshape Proof. }}%
  {{\slshape #1.} }%
 }%
{\qed\bigskip}

\newcounter{alphabet}
\newcounter{tmp}
\newenvironment{Thm}[1][]{\refstepcounter{alphabet}%
\bigskip%
\noindent%
{\bf Theorem \Alph{alphabet}}%
\ifthenelse{\equal{#1}{}}{}{ (#1)}%
{\bf .} \itshape}{\vskip 8pt}
\makeatletter
\newcommand{\Ref}[1]{\@ifundefined{r@#1}{}{\setcounter{tmp}{\ref{#1}}\Alph{tmp}}}
\makeatother
\newenvironment{Lem}[1][]{\refstepcounter{alphabet}%
\bigskip%
\noindent%
{\bf Lemma \Alph{alphabet}}%
{\bf .} \itshape}{\vskip 8pt}

\newcommand{\ID}{{\mathbb D}}

\def\be{\begin{equation}}
\def\ee{\end{equation}}

\newcommand{\bee}{\begin{enumerate}}
\newcommand{\eee}{\end{enumerate}}

\newcommand{\blem}{\begin{lem}}
\newcommand{\elem}{\end{lem}}
\newcommand{\bthm}{\begin{thm}}
\newcommand{\ethm}{\end{thm}}
\newcommand{\bcor}{\begin{cor}}
\newcommand{\ecor}{\end{cor}}
\newcommand{\beg}{\begin{example}}
\newcommand{\eeg}{\end{example}}
\newcommand{\begs}{\begin{examples}}
\newcommand{\eegs}{\end{examples}}
\newcommand{\bdefe}{\begin{defn}}
\newcommand{\edefe}{\end{defn}}
\newcommand{\bprob}{\begin{prob}}
\newcommand{\eprob}{\end{prob}}
\newcommand{\bques}{\begin{ques}}
\newcommand{\eques}{\end{ques}}
\newcommand{\bei}{\begin{itemize}}
\newcommand{\eei}{\end{itemize}}
\newcommand{\bcon}{\begin{conj}}
\newcommand{\econ}{\end{conj}}
\newcommand{\bcons}{\begin{conjs}}
\newcommand{\econs}{\end{conjs}}
\newcommand{\bprop}{\begin{propo}}
\newcommand{\eprop}{\end{propo}}
\newcommand{\br}{\begin{rem}}
\newcommand{\er}{\end{rem}}
\newcommand{\brs}{\begin{rems}}
\newcommand{\ers}{\end{rems}}
\newcommand{\bo}{\begin{obser}}
\newcommand{\eo}{\end{obser}}
\newcommand{\bos}{\begin{obsers}}
\newcommand{\eos}{\end{obsers}}
\newcommand{\bpf}{\begin{pf}}
\newcommand{\epf}{\end{pf}}
\newcommand{\ba}{\begin{array}}
\newcommand{\ea}{\end{array}}
\newcommand{\beq}{\begin{eqnarray}}
\newcommand{\beqq}{\begin{eqnarray*}}
\newcommand{\eeq}{\end{eqnarray}}
\newcommand{\eeqq}{\end{eqnarray*}}

\newcommand{\ds}{\displaystyle}

\newcounter{minutes}
\divide\time by 60
\newcounter{hours}
\multiply\time by 60 \addtocounter{minutes}{-\time}
\begin{document}

\bibliographystyle{amsplain}

%
%\begin{center}
%{\tiny \texttt{FILE:~\jobname .tex,
%        printed: \number\year-\number\month-\number\day,
%        \thehours.\ifnum\theminutes<10{0}\fi\theminutes}
%}
%\end{center}

\title[Coefficients of bounded functions]{New Inequalities for the Coefficients of unimodular bounded Functions}
% \title[Logarithmic coefficients inequalities]{Logarithmic coefficients inequalities
%for certain subfamilies of univalent functions}

%=========================================================================
\thanks{%$^\dagger$
File:~\jobname .tex,
          printed: \number\day-\number\month-\number\year,
          \thehours.\ifnum\theminutes<10{0}\fi\theminutes}
%=========================================================================

\author[S. Ponnusamy]{Saminathan Ponnusamy
%$^\dagger $
%${}^{~\mathbf{*}}$
}
\address{
S. Ponnusamy, Department of Mathematics,
Indian Institute of Technology Madras, Chennai-600 036, India.
}
\email{samy@iitm.ac.in}

\author[R. Vijayakumar]{Ramakrishnan Vijayakumar}
\address{
R. Vijayakumar, Department of Mathematics,
Indian Institute of Technology Madras, Chennai-600 036, India.
}
\email{mathesvijay8@gmail.com}

\author[K.-J. Wirths]{Karl-Joachim Wirths}
\address{K.-J. Wirths, Institut f\"ur Analysis und Algebra, TU Braunschweig,
38106 Braunschweig, Germany.}
\email{kjwirths@tu-bs.de}

\keywords{Analytic functions, Schwarz lemma, Bohr's inequality, unimodular bounded functions
%subordination and quasisubordination
%\\
%$%{}^{\mathbf{*}}
%^\dagger$ {\tt This author is on leave from the Department of Mathematics,
%Indian Institute of Technology Madras, Chennai-600 036, India}
}
%\thanks{The authors thank the referees for their careful reading of the paper.
%The first author is on leave from IIT Madras, and is currently at ISI Chennai.
%}
%\date{\today  %June. 30, 09
%;  File: 2013(2).tex}

\begin{abstract}

The classical inequality of Bohr asserts that if a power series converges in the unit disk and its sum has modulus
less than or equal to $1$, then the sum of absolute values of its terms is less than or equal to $1$ for the subdisk
$|z|<1/3$ and $1/3$ is the best possible constant. Recently, there has been a number of investigations on this topic.
In this article, we present related inequalities using $\sum_{n=0}^{\infty}|a_n|^2r^{2n}$ that generalize for example
the well known inequality
\[\sum_{n=0}^{\infty}|a_n|^2r^{2n} \leq 1.
\]
\end{abstract}

\maketitle
\pagestyle{myheadings}
\markboth{S. Ponnusamy, R. Vijayakumar and K.-J. Wirths}{Coefficients of bounded functions}

\section{Introduction and Main results}

Let $\mathbb{D}=\{z\in \mathbb{C}:\,|z|<1\}$ denote the open unit disk and $\mathcal B$ denote the
class of all analytic functions in $\ID$ such that $|f(z)|\leq 1$ in $\ID$.  Let us first recall the theorem of H. Bohr \cite{Bohr-14} in 1914,
which inspired a lot in the recent years. Bohr originally discovered this inequality with $r\leq 1/6$,
and Bohr's theorem was proved in this form by  F. Wiener, and later
the same was obtained independently by M. Riesz, I. Schur and some others. See for example \cite{AAPon1,GarMasRoss-2018}.

\begin{Thm}\label{PVW1-theA}
If $f\in {\mathcal B}$ and $f(z)=\sum_{n=0}^{\infty} a_n z^n$,
then
\be\label{KKP1-eq1}
\sum_{n=1}^{\infty} |a_n|\, r^n \leq 1 -|a_0| %~(=\dist (f(0),\partial \ID))
~\mbox{ for  $r\leq 1/3$}
\ee
and  the constant $1/3$ cannot be improved.
\end{Thm}

The inequality \eqref{KKP1-eq1} is known as the classical Bohr inequality and the number $1/3$ is called the Bohr radius
for the family ${\mathcal B}$. It is worth pointing out that if $|a_0|$ in \eqref{KKP1-eq1} is replaced by $|a_0|^2$, then the constant $1/3$
could be replaced by $1/2$.
Moreover, if $a_0=0$ in Theorem \Ref{PVW1-theA} then the sharp Bohr radius is known to be $1/\sqrt{2}$ (see for example \cite{Bom-62}, \cite{KayPon1}
and \cite[Corollary 2.9]{PaulPopeSingh-02-10}).
Extensions, modifications and improvements of this result can be found in many recent papers
\cite{AliBarSoly,BDK5,BhowDas-18,EvPoRa-2017, KayPon1,KayPon2,KayPonShak1,LP2018,LLP2020,LSX2018,PaulSingh-04-11,PaulSingh-06-12}.
For example,  several improved versions of Theorem \Ref{PVW1-theA} are established recently in \cite{AlkKayPon1,KayPon3}.
The true outburst of activity around Bohr's estimate occurred after publication of the article of
Boas and Khavinson  \cite{BoasKhavin-97-4}
and Dixon \cite{Dixon-95-7} used Bohr's original theorem
to construct a Banach algebra which is not an operator algebra, yet satisfies the non-unital von Neumann's inequality.
Paulsen et al. \cite{PaulPopeSingh-02-10}
have applied operator-theoretic techniques to extend it to Banach algebras and obtain  multidimensional generalizations of Bohr's inequality.
Further investigations such as interconnections among multidimensional Bohr radii and local Banach space theory may be obtained from the
survey articles of Abu-Muhanna et al. \cite{AAPon1}, Defant and Prengel \cite{DePre06}, Garcia et al. \cite{GarMasRoss-2018}.
%See also  Defant et al. \cite{DeGarM04} and the references therein.

The proof of the inequality  \eqref{KKP1-eq1} used the   sharp coefficient inequalities which may be obtained as an
application of Pick's invariant form of Schwarz's lemma for $f\in {\mathcal B}$:
$$|f'(z)|\leq \frac{1-|f(z)|^2}{1-|z|^2}, \quad z\in\ID.
$$
In particular, $|f'(0)|=|a_1|\leq 1-|f(0)|^2=1-|a_0|^2$ from which  sharp inequalities $|a_n| \leq 1-|a_0|^2$
$(n\geq 1,~ f\in {\mathcal B})$ follow.
As remarked in \cite{KayPon1}, we were not able to obtain sharp result due to the fact that in the extremal case  $|a_0|<1$. Lately in
\cite{AlkKayPon1}, a sharp version of Theorem \Ref{PVW1-theA} has been achieved and this works for any individual function
from ${\mathcal B}$.

Our aim is the comparison of $\sum_{n=0}^{\infty}|a_n|r^n$ with another functional often considered in function theory, namely
$\sum_{n=0}^{\infty}|a_n|^2r^{2n}$ that will be abbreviated as $\|f\|_r^2$. Our first result is a generalization of the well known inequality
$  \|f\|_r^2\,\leq 1$ for $f\in \mathcal{B}$.

\begin{thm}\label{PVW1-Th1}
Suppose that $f\in {\mathcal B}$ and  $f(z)=\sum_{n=0}^{\infty} a_n z^n$. Then the inequality
\[ \sum_{n=0}^{\infty}|a_n|r^n\,\leq \,\frac{1}{1-r}\left(1 -r\|f\|_r^2\right)
\]
is valid for $r\in [0,1)$. Equality is attained for $f(z)=1$, $z\in \ID.$
\end{thm}

Further we prove the following results that can be looked as refinements of Bohr's inequality.

\begin{thm}\label{PVW1-Lem1}
Suppose that $f\in {\mathcal B}$, $f(z)=\sum_{n=0}^{\infty} a_n z^n$,
and $f_0(z)=f(z)-a_0$.  Then
\begin{equation}\label{PVW1-eq1}
 \sum_{n=0}^\infty |a_n|r^n+  \left(\frac{1}{1+|a_0|}+\frac{r}{1-r}\right)\|f_0\|_r^2 \leq 1 ~\mbox{ for  }~ r \leq \frac{1}{2+|a_0|}
\end{equation}
and the numbers $\frac{1}{2+|a_0|}$ and $\frac{1}{1+|a_0|}$ cannot be improved.  Moreover,
\begin{equation}\label{PVW1-eq2}
|a_0|^2+ \sum_{n=1}^\infty |a_n|r^n+ \left(\frac{1}{1+|a_0|}+\frac{r}{1-r}\right)\|f_0\|_r^2 \leq 1 ~\mbox{ for  }~ r \leq \frac{1}{2}
\end{equation}
and the numbers $\frac{1}{2}$ and $\frac{1}{1+|a_0|}$ cannot be improved.
\end{thm}

We remark that $\frac{1}{3}\leq \frac{1}{2+|a_0|}\leq \frac{1}{2}$. In  particular, in the case of $a_0=0$ the conclusion of  Theorem \ref{PVW1-Lem1}, namely, the inequality \eqref{PVW1-eq1} gives that
$$ \sum_{n=1}^\infty |a_n|r^n+  \frac{1}{1-r}\|f\|_r^2 \leq 1 ~\mbox{ for  }~ r \leq \frac{1}{2}.
$$
Unfortunately, the direct substitution for $a_0$ would not yield sharp Bohr radii.
In the following, we prove two different improved versions of the last inequality.

\bthm\label{PVW1-Lem5}
Suppose that $f\in {\mathcal B}$ and $f(z)=\sum_{n=1}^{\infty} a_n z^n$. Then we have the following:
\begin{enumerate}
\item[{\rm (a)}] $\ds \sum_{n=1}^{\infty}|a_n|r^n + \left(\dfrac{1}{1+|a_1|}+\dfrac{r}{1-r}\right)\sum_{n=2}^{\infty}|a_{n}|^2r^{2n-1}
\leq  1 ~\mbox{ for }~ r\leq \dfrac{3}{5}$.\\ The number $3/5$ is sharp.
\item[{\rm (b)}] $\ds \sum_{n=1}^{\infty}|a_n|r^n +   \left(\dfrac{r^{-1}}{1+|a_1|}+\dfrac{1}{1-r}\right)\|f\|_r^2 \leq  1 ~\mbox{ for }~ r\leq~ \frac{5-\sqrt{17}}{2}.$\\ The number $\frac{5-\sqrt{17}}{2}$ is sharp.

%\item[{\rm (c)}] Let  $|a_1|\,=\,\frac{1}{\sqrt{2}}$. Then
%    $\ds \sum_{n=1}^{\infty}|a_n|r^n +   \left(\dfrac{r^{-1}}{1+|a_1|}+\dfrac{1}{1-r}\right)\|f\|_r^2 \leq  1 ~\mbox{ for }~ r\leq~ \frac{1}{2}.$
%    \\The number $\frac{1}{2}$ is sharp.
 \item[{\rm (c)}]
     $\ds \sum_{n=1}^{\infty}|a_n|r^n +   \left(\dfrac{r^{-1}}{1+|a_1|}+\dfrac{1}{1-r}\right)\|f\|_r^2 \leq  1 ~\mbox{ for $r \leq r(a)$,}$ where
$$
 r(a) = \dfrac{2(1+a)}{1+2a+2a^2\,+\, \sqrt{4a^4+8a+5}} .
$$
The radius $r(a)$ is sharp for any $a\in [0,1)$.
\end{enumerate}	
\ethm

We remark that
$$\frac{3}{5}>\frac{1}{2}=r\left (\frac{1}{\sqrt{2}}\right )>\frac{5-\sqrt{17}}{2}.
$$
The proofs of these results will be presented in Section \ref{Sec2-PVW2}.

\section{Proofs of the Theorems}\label{Sec2-PVW2}

The following lemma due to Carlson \cite{Carlson-40} is key for the proofs of Theorems \ref{PVW1-Th1}, \ref{PVW1-Lem1} and \ref{PVW1-Lem5}.

\begin{Lem}\label{PVW1-lem2}
%{\rm \cite{Carlson-40}}
Suppose that $f\in {\mathcal B}$ and $f(z)=\sum_{n=0}^{\infty} a_n z^n$.
%Let $f(z)=\sum_{n=0}^{\infty}a_{n}z^n$ be an analytic f and $|f(z)|\leq  1$ in $\mathbb{D}.$
Then the following inequalities hold.
\begin{enumerate}
	\item[{\rm (a)}] $|a_{2n+1}|\leq 1-|a_0|^2-\cdots - |a_n|^2,\ n=0,1,\ldots$ %\label{hyp1}
	\item[{\rm (b)}] $|a_{2n}|\leq 1-|a_0|^2-\cdots -|a_{n-1}|^2 - \frac{|a_n|^2}{1+|a_0|} ,\ n=1,2,\ldots$. %\label{hyp2}
\end{enumerate}
Further, to have equality in {\rm (a)} it is necessary that $f$ is a rational function of the form
%$$ f(z)=\frac{a_{0}+a_{1}z+ \cdots + a_{n}z^{n}+z^{2n+1}}{1+\overline{a_n}z^n+ \cdots +\overline{a_0}z^{2n+1}}
%$$
 $$ f(z)=\frac{a_{0}+a_{1}z+ \cdots + a_{n}z^{n}+\epsilon z^{2n+1}}{1+(\overline{a_n}z^n+ \cdots +\overline{a_0}z^{2n+1})\epsilon},\quad |\epsilon|=1,
$$
and  to have equality in {\rm (b)} it is necessary that $f$ is a rational function of the form
$$ f(z)=\frac{a_{0}+a_{1}z+ \cdots + \frac{a_{n}}{1+|a_0|} z^{n}+ \epsilon z^{2n} }{1+\left(\frac{\overline{a_n}}{1+|a_0|}z^n+ \cdots +\overline{a_0}z^{2n}\right) \epsilon},
\quad |\epsilon|=1,
$$
where the condition $a_0 \overline{a_n}^2 \epsilon$ is non-positive real.
\end{Lem}

\subsection{Proof of Theorem \ref{PVW1-Th1}}
Let $M_f(r)=\sum_{n=0}^{\infty}|a_n|r^n$ denote the  associated majorant series for $f(z)=\sum_{n=0}^{\infty}a_nz^n$ which is analytic  and $|f(z)| \leq 1$ in $\ID$.
Using Lemma \Ref{PVW1-lem2}, one has
\beq\label{eq-th1a}
M_f(r)&= &|a_0|+\sum_{n=1}^{\infty}|a_{2n}|r^{2n}+\sum_{n=0}^{\infty}|a_{2n+1}|r^{2n+1}\nonumber\\
&\leq&  |a_0|+\sum_{n=1}^{\infty}\left[1-\sum_{k=0}^{n-1}|a_k|^2-\frac{|a_n|^2}{1+|a_0|}\right]r^{2n}
+\sum_{n=0}^{\infty}\left[1-\sum_{k=0}^{n}|a_k|^2\right]r^{2n+1}\nonumber\\ &=&|a_0|+\sum_{n=1}^{\infty}r^n-|a_0|^2\left(\sum_{n=1}^{\infty}r^n\right)-\left(\frac{1}{1+|a_0|}
+\sum_{n=1}^{\infty}r^n\right)\sum_{n=1}^{\infty}|a_n|^2r^{2n}\nonumber\\
&=&|a_0|+ \frac{r}{1-r} \left(1-|a_0|^2\right)-\left(\frac{1}{1+|a_0|}+\frac{r}{1-r}\right) \|f_0\|_r^2,
\eeq
 where $f_0(z)= f(z)-f(0)$.
Now, we consider the case $a_0=0.$ Inserting this into \eqref{eq-th1a} delivers
\[\sum_{n=1}^{\infty}|a_n|r^n \leq \frac{1}{1-r}\left(r - \|f_0\|_r^2\right).
\]
Due to the  classical lemma of Schwarz, we have $f(z) = z\omega(z)$, $\omega \in {\mathcal B}.$ If we let
\[\omega(z) = \sum_{n=0}^{\infty}\alpha_nz^n,
\]
then we have $a_{n+1} = \alpha_n$ for $n\geq 0$. Inserting this equation into the above inequality and dividing by $r$ implies
\[M_\omega(r) = \sum_{n=0}^{\infty}|\alpha_n|r^n \leq \frac{1}{1-r}\left(1 - r\sum_{n=0}^{\infty}|\alpha_n|^2r^{2n}\right) =
\frac{1}{1-r}\left(1 - r\|\omega \|_r^2\right)
\]
for any $\omega \in {\mathcal B}$, $r\in [0,1).$ The assertion on equality is trivial. This completes the proof of Theorem \ref{PVW1-Th1}.
\hfill $\Box$

\bigskip
%  \noindent

\subsection{Proof of Theorem \ref{PVW1-Lem1}}
 Since $r/(1-r)$ is an increasing function of $r\in (0,1)$, it follows from  the relation \eqref{eq-th1a} that for $r \leq 1/(2+|a_0|)$
\beqq
 M_f(r)+\left(\frac{1}{1+|a_0|}+\frac{r}{1-r}\right) \|f_0\|_r^2 &\leq & |a_0|+ \frac{r}{1-r} \left(1-|a_0|^2\right) =:\Psi (r)\\
 &\leq &\Psi \left ( \frac{1}{2+|a_0|}\right )=1
\eeqq
which implies the inequality \eqref{PVW1-eq1}.

To prove that the radius is sharp, we consider the function $f=\varphi_a$ given by
$$\varphi_a (z)  = \frac{a-z}{1-az} =a - (1-a^2)\sum_{k=1}^\infty a^{k-1} z^{k}, \quad z\in\ID,
$$
where $a\in (0,1)$. For this function, with $a_0=a$ and $a_n=-(1-a^2)a^{n-1}$, straightforward calculations show that
$$M_{\varphi_a }(r)=\sum_{n=0}^\infty |a_n|r^n =a+r\frac{1-a^2}{1-ar}
$$
and
$$\left(\frac{1}{1+a}+\frac{r}{1-r}\right)\|\varphi_a -a\|_r^2= \frac{1+ar}{(1-r)(1+a)}\frac{(1-a^2)^2r^2}{1-a^2r^2}
= \frac{(1-a)^2(1+a)r^2}{(1-r)(1-ar)}
$$
so that

\vspace{8pt}

$\ds M_{\varphi_a }(r)+ \left(\frac{1}{1+a}+\frac{r}{1-r}\right)\|\varphi_a -a\|_r^2
$

\beq\label{eq-th1b}
%\sum_{k=0}^\infty |a_k|r^k +\left(\frac{1}{1+|a_0|}+\frac{r}{1-r}\right)\sum_{n=1}^{\infty}|a_n|^2r^{2n}
%M_{\varphi_a }(r)+ \left(\frac{1}{1+|a_0|}+\frac{r}{1-r}\right)\|\varphi_a -a\|_r^2
&=& a+r\frac{1-a^2}{1-ar}+  \frac{(1-a)^2(1+a)r^2}{(1-r)(1-ar)}\nonumber\\
&=&1-(1- a)\left [\frac{1-(1+2a)r}{1-ar} - \frac{(1-a^2)r^2}{(1-r)(1-ar)}\right ]\nonumber\\
&=&1+\frac{(1- a)[(2+a)r-1]}{1-r}
\eeq
which shows that the left hand side is bigger than $1$ whenever $r>1/(2+a)$.

The second part, namely the inequality \eqref{PVW1-eq2}, is clear if we replace $|a_0|$ by $|a_0|^2$
in the majorant sum, and thus, from \eqref{eq-th1a}, we see that
$$|a_0|^2+ \sum_{n=1}^\infty |a_n|r^n+ \left(\frac{1}{1+|a_0|}+\frac{r}{1-r}\right)\|f_0\|_r^2\leq |a_0|^2 +(1-|a_0|^2)\frac{r}{1-r}
$$
which is obviously less than or equal to $1$ whenever $r\leq 1/2$. The sharpness of the constant $1/2$ can be established as
in the previous case. Indeed, for the function $f=\varphi _a$, we have by \eqref{eq-th1b} that
\beqq
M_{\varphi_a }(r)+  \left(\frac{1}{1+a}+\frac{r}{1-r}\right)\|\varphi_a -a\|_r^2-a+a^2
&=& -a(1-a) + 1+\frac{(1- a)[(2+a)r-1]}{1-r}\\
&=&1 + \frac{(1- a^2)(2r-1)}{1-r}
\eeqq
which is bigger than $1$ whenever $r>1/2$. The proof of the theorem is complete.
\hfill $\Box$

%\epf

\subsection{Proof of Theorem \ref{PVW1-Lem5}}

Let $f(z)=\sum_{n= 1}^\infty a_{n}z^{n},$ where $|f(z)|\leq 1$ for  $z\in\ID.$ At first, we remark that the function $f$ can be represented as
$f(z)=z g(z),$ where  $g(z)=\sum_{n = 0}^\infty b_{n}z^{n}$ is analytic in $\ID$ with $b_n = a_{n+1}$ and $|g(z)|\leq 1$ for  $z\in\ID$.
Let $|b_0|=|a_1|=a$ and $g_0(z)=g(z)-b_0$.  Then it follows from the proof of Theorem \ref{PVW1-Lem1} that
\be\label{PVW1-eq4a}
\sum_{n=1}^{\infty}|a_n|r^n   =  \sum_{n=0}^{\infty}|b_n|r^{n+1}
 \leq  r \left[ a+ \dfrac{r}{1-r} \left(1-a^2\right)-\left(\dfrac{1}{1+a}+\dfrac{r}{1-r}\right)\|g_0\|_r^2 \right]
\ee
and thus, we have
\be\label{PVW1-eq4}
\sum_{n=1}^{\infty}|a_n|r^n + \left(\dfrac{1}{1+a}+\dfrac{r}{1-r}\right)\sum_{n=1}^{\infty}|b_n|^2r^{2n+1}  \leq  \Psi(a;r) ,
%r \left[ a+\dfrac{r}{1-r}\left(1-a^2\right) \right]
\ee
where
$$ \Psi(x;r) = rx+\dfrac{r^2}{1-r}\left(1-x^2\right) \mbox{ for $x \in [0,1]$ and $r \in [0,1)$} .
$$
We just need to maximize $\Psi(x;r)$  (with respect to $x$) over the interval $[0,1].$
We see that $\Psi$ has a critical point at
$x_0=(1-r)/(2r)$ and obtain that the maximum occurs at this point so that
\be\label{PVW1-eq5}
\Psi (x;r) \leq \Psi ( x_0;r)= \dfrac{1-r}{2}+\dfrac{(3r-1)(1+r)}{4(1-r)}=1- \dfrac{(3-5r)(1+r)}{4(1-r)},
\ee
which is less than or equal to $1$ whenever $r\leq 3/5$. The first inequality {\rm (a)} follows from \eqref{PVW1-eq4} and \eqref{PVW1-eq5}.

To discuss the sharpness in \text{(a)}, we consider the function $f=\varphi_a$ given by
$$\varphi_a (z)  = z\left(\frac{a-z}{1-az}\right) =az - (1-a^2)\sum_{n=1}^\infty a^{n-1} z^{n+1}, \quad z\in\ID,
$$
where $a\in (0,1)$. For this function, with $a_1=a$ and $a_n=-(1-a^2)a^{n-2}$ for $n\geq 2$, straightforward calculations show that
$$M_{\varphi_a }(r)=\sum_{n=1}^\infty |a_n|r^n =ar+(1-a^2)  \frac{r^{2}}{1-ar}
$$
and
$$\left(\dfrac{1}{1+a}+\dfrac{r}{1-r}\right)\sum_{n=2}^\infty |a_n|^{2}r^{2n-1}= \frac{(1-a^2)(1-a)r^3}{(1-r)(1-ar)}
$$
so that the sum gives
\beqq
%\sum_{k=0}^\infty |a_k|r^k +\left(\dfrac{1}{1+|a_0|}+\dfrac{r}{1-r}\right)\sum_{n=1}^{\infty}|a_n|^2r^{2n}
M_{\varphi_a }(r)+ \left(\dfrac{1}{1+a}+\dfrac{r}{1-r}\right)\sum_{n=2}^\infty |a_n|^{2}r^{2n-1}
%&=& ar+r^{2}\left(\frac{1-a^2}{1-ar}\right)+\frac{(1-a^2)(1-a)r^3}{(1-r)(1-ar)}\\
%&=&ar+r^{2}\left(\frac{1-a^2}{1-ar}\right)\left[1+\frac{r(1-a)}{1-r}\right]\\
&=&ar+ (1-a^2)\frac{r^{2}}{1-r} =\frac{\Phi(a,r)}{1-r},
\eeqq
where $ \Phi(a,r)=ar(1-r)+(1-a^2)r^2$. We want to show that $ \Phi(a,r)>1-r$ for any $r>3/5$ and an appropriate
number $a\in (0,1).$ Since %the partial derivative of $\Phi$ with respect to $a$  equals
\[ \frac{\partial \Phi(a,r)}{\partial a}= r(1-r) - 2ar^2,
\]
we see that $\Phi$ is an increasing function of $a$ for $a \leq (1-r)/(2r)=:c(r)$.
For $r \in (1/3, 1), $  we have $c(r) \in (0,1).$ Now, we calculate
\[ \Phi\left(\frac{1-r}{2r},r\right) = \frac{(1-r)^2 +4r^2}{4},
\]
and we see that this quantity is bigger than $1-r$ (which is equivalent to $5r^2+2r-3=(5r-3)(r+1)>0$) for $r>3/5.$
This completes the proof of the sharpness of the constant $3/5$ in \text{(a)}.

Next, we prove the inequalities (b) and (c). In order to do this, we use the abbreviation $|a_1|\,=\,a$, and we  make minor changes in the proof of (a).
 Accordingly, from \eqref{PVW1-eq4}, we find that

\vspace{8pt}

$\ds \sum_{n=1}^{\infty}|a_n|r^n + \dfrac{1}{1-r} \|f\|_r^2
$
\begin{align*}
&\leq
% r \left[ a+ \dfrac{r}{1-r} \left(1-a^2\right)-\left(\dfrac{1}{1+a}+\dfrac{r}{1-r}\right)\|g_0\|_r^2 \right]
%+ \dfrac{1}{1-r}  \sum_{n=1}^{\infty}|a_n|^2r^{2n} \\
%&=
r \left[ a+ \dfrac{r}{1-r} \left(1-a^2\right)-\left(\dfrac{1}{1+a}+\dfrac{r}{1-r}\right)\sum_{n=2}^{\infty}|a_n|^2r^{2(n-1)} \right]
 + \dfrac{r^2}{1-r} \sum_{n=1}^{\infty}|a_n|^2r^{2(n-1)}\\
&= ra + \dfrac{r^2}{1-r}  -\dfrac{r^{-1}}{1+a} \sum_{n=2}^{\infty}|a_n|^2r^{2n} \\
&= ra +  \dfrac{r^2}{1-r} + \dfrac{r  a^2}{1+a}   - \dfrac{r^{-1}}{1+a} \sum_{n =1}^{\infty}|a_n|^2r^{2n}
 \end{align*}
and thus, we have \be\label{PVW1-eq6}
 \sum_{n=1}^{\infty}|a_n|r^n +  \left (\dfrac{r^{-1}}{1+|a_1|}+\dfrac{1}{1-r}\right )\|f\|_r^2 \leq \Xi(a;r) ,
 \ee
where
\[
\Xi(a;r) =ra +\dfrac{r^2}{1-r}+\dfrac{r  a^2}{1+a} = 1- \dfrac{B(a,r)}{(1-r)(1+a)},
\]
where $B(a,r)=r^2 (2a^2 -1)-r(1+ 2a+ 2a^2 ) + 1+a.$ Since
\[
\frac{\partial \Xi(a,r)}{\partial a} = r + \frac{2ar+a^2r}{(1+a)^2},
\]
the function $\Xi$ is a strictly monotonic increasing function of $a$. Therefore
\[
\Xi(a,r) \leq \Xi(1,r) = \frac{3r - r^2}{2(1 - r)} \leq 1,
\]
if $2 - 5r + r^2 \geq 0$. This is the case for $r \in [0, (5-\sqrt{17})/2]$. To see that this radius is sharp,
we consider the function $f(z)=z$. For this function the left hand side of  (b) takes the value
\[\frac{3r - r^2}{2(1 - r)}
\]
which is less than or equal to unity if and only if $r$ takes the values indicated above. This completes the
proof of part (b).

To prove (c) we consider
\[  \frac{\partial \Xi(a,r)}{\partial r} = a + \frac{r(2-r)}{(1-r)^2} +\frac{a^2}{1+a} >0.
\]
Therefore $\Xi(a,r)$ is a strictly monotonic increasing function of $r$.
Hence, for a fixed value of $a$, the desired inequalities are valid for $r\in [0,r(a)]$, where
$r(a)$ is the first positive zero of the function $B(a,r)$.
%\[B(a,r)=r^2 (2a^2 -1)-r(1+ 2a+ 2a^2 ) + 1+a.
%\]
If $a=1/\sqrt{2}$, the only zero of the function $B$ is $r=1/2$.
%This proves the inequality (c).

For $a\in (1/\sqrt{2},1)$, the solution $r_{\pm}$ of the equation $B(a,r) = 0$ given by
\beqq
 r_{\pm} &=& \dfrac{(1+2a+2a^2)\pm \sqrt{(1+2a+2a^2)^2 - 4(2a^2-1)(1+a)}}{2(2a^2-1)} \\
 &=& \dfrac{(1+2a+2a^2)\pm \sqrt{4a^4+8a+5}}{2(2a^2-1)}
\eeqq
are both positive. Hence, we can take
$$r(a) = r_{-}=\dfrac{(1+2a+2a^2)\, - \,\sqrt{4a^4+8a+5}}{2(2a^2-1)}= \dfrac{2(1+a)}{1+2a+2a^2\,+\, \sqrt{4a^4+8a+5}} \,.
$$
The last expression is convenient to allow $a=1/\sqrt{2}$ as well, since the expression on the right
gives $1/2$ when $a=1/\sqrt{2}$, as remarked above.

Finally, for $a \in (0,1/\sqrt{2})$, we have
\[
\left.\frac{\partial B(a,r)}{\partial r}\right|_{r=0} < 0.
\]
This proves that one of the zeros is negative. The other one is again $r_{-}$ such that we have $r(a) = r_{-}$.
This completes the proof of the inequality in (c).

Now, it remains to prove the sharpness in the assertions (c). To that end, we again consider the
function $f=\varphi_a$ from part (a) for the values $a \in [0,1)$. If we add $M_{\varphi_a }(r)$ and
$$\left(\frac{r^{-1}}{1+a}+\frac{1}{1-r}\right)\sum_{n=1}^\infty |a_n|^{2}r^{2n}
= \left(\frac{1+ar}{r(1+a)(1-r)}\right)\left(a^2r^2+\frac{(1-a^2)^2r^4}{1-a^2r^2}\right),
$$
then we get

\vspace{8pt}
$\ds M_{\varphi_a }(r)+ \left(\frac{r^{-1}}{1+a}+\frac{1}{1-r}\right)\sum_{n=1}^\infty |a_n|^{2}r^{2n}$
\beqq
%\sum_{k=0}^\infty |a_k|r^k +\left(\dfrac{1}{1+|a_0|}+\dfrac{r}{1-r}\right)\sum_{n=1}^{\infty}|a_n|^2r^{2n}
&=&ar+r^{2}\left(\frac{1-a^2}{1-ar}\right)+\left(\frac{1+ar}{r(1+a)(1-r)}\right)\left(a^2r^2+\frac{(1-a^2)^2r^4}{1-a^2r^2}\right)\\
&=&ar+r^{2}\left(\frac{1-a^2}{1-ar}\right)\left(1+\frac{r(1-a)}{1-r}\right)+a^2r \left(\frac{1+ar}{(1+a)(1-r)}\right)\\
&=&ar+\frac{r^2(1-a^2)}{1-r}+a^2r \left(\frac{1+ar}{(1+a)(1-r)}\right)
\eeqq
which is greater than 1 whenever $B(a,r) < 0.$

Therefore, the above discussion of $B$ as a function of $r$ for a fixed value of $a$ immediately delivers that the radius
in (c) is sharp. This completes the proof of Theorem \ref{PVW1-Lem5}.
\hfill $\Box$

%\epf

\br
It may be interesting that for $a_1=0$ the sharp radius in the second part of Theorem 2 has the value
$(\sqrt{5}-1)/2$ and that in this case the extremal function is given by $\varphi_a(z)=z^2.$
\er

\subsection*{Acknowledgments}
%The authors thank the referee for his/her useful suggestions for the improvement of the presentation.
The work of the first author is supported by Mathematical Research Impact Centric Support of DST, India  (MTR/2017/000367).
The authors want to thank H. Schellwat from \"Orebro University for his help in getting the article of
F.  Carlson \cite{Carlson-40}.

\end{document}